\documentclass[twocolumn, 11pt]{article}
\RequirePackage[OT1]{fontenc}
\RequirePackage{amsthm,amsmath}
\RequirePackage[numbers]{natbib}

\newcommand{\rr}{{\Bbb R}}

\newcommand{\one}{{\hbox{1{\kern -0.35em}1}}}

\newcommand{\M}{{\cal M}}

\newcommand{\beq}[1]{\begin{eqnarray} \label{#1}}
\newcommand{\eeq}{\end{eqnarray}}
\newcommand{\bed}{\begin{displaymath}}
\newcommand{\eed}{\end{displaymath}}
\newcommand{\bea}{\bed\begin{array}{rl}}
\newcommand{\eea}{\end{array}\eed}
\newcommand{\disp}{\displaystyle}
\newcommand{\al}{\alpha}

\newtheorem{thm}{Theorem}[section]

\newtheorem{lem}[thm]{Lemma}
\newtheorem{rem}[thm]{Remark}

\newtheorem{defn}[thm]{Definition}

\newcommand{\thmref}[1]{Theorem~{\rm \ref{#1}}}
\newcommand{\lemref}[1]{Lemma~{\rm \ref{#1}}}

\def\openbox{$\sqcup\llap{$\sqcap$}$}
\def\endproof{\unskip \enskip
    \null \nobreak \hfill \openbox \par}

\newcommand{\ad}{&\!\!\!\disp}
\newcommand{\aad}{&\disp}
\newcommand{\barray}{\begin{array}{ll}}
\newcommand{\earray}{\end{array}}

\usepackage{makeidx}  
\usepackage{amsfonts}
\makeindex            
\usepackage{times}  
\usepackage{graphicx} 

\begin{document}
 


\date{}

\title{Optimal Oil Production and Taxation in Presence of Global Disruptions}
\author{Moustapha Pemy \thanks{Department of Mathematics, Towson
University, Towson, MD 21252-0001, mpemy@towson.edu
 }}
\maketitle


\begin{abstract} 
This paper studies the optimal extraction policy of an oil field as well as the efficient taxation of the revenues generated. Taking into account the fact that the oil price in worldwide commodity markets fluctuates randomly following global and seasonal macroeconomic parameters, we model the evolution of the oil price as a mean reverting regime-switching jump diffusion process. Given that oil producing countries rely on oil sale revenues as well as taxes levied on oil companies for a good portion of the revenue side of their budgets, we formulate this problem as a differential game where the two players are the mining company whose aim is to maximize the revenues generated from its extracting activities and the government agency in charge of regulating and taxing natural resources. We prove the existence of a Nash equilibrium and the convergence of an approximating scheme for the value functions. Furthermore, optimal extraction and fiscal policies that should be applied when the equilibrium is reached are derived.
 A numerical example is presented to illustrate these results.
\vskip 0.2 in

\end{abstract}





\section{Introduction}
Oil and natural gas have always been the main sources of revenues for a large number of developing countries as well as some industrialized countries around the world. Oil extraction policies vary from a country to another. In some countries, the extraction is done by a state-owned corporation in others it is done by foreign multinationals. 
 The earliest work on the extraction of natural resources was done by Hotelling (1931) who derived an optimal extraction policy under the assumption that the commodity price is constant. 
Many economists have proposed various extensions of the Hotelling model by taking into account the uncertainty and randomness of commodity prices. 

 The majority of oil extraction contracts signed between multinational oil companies and governments of oil-rich nations are in the form of profit sharing agreements where each party will take a fraction of the profits.  In addition, the host nation is also entitled to collect tax from all companies. This creates a very interesting dynamic for these two parties with converging as well as conflicting interests during the lifetime of the mining contract.   
 We formulate this problem as a differential game where the two players are the multinational oil company and the government. 
 To the best of our knowledge, this is the first time this approach is used to characterize the interplay between extraction and taxation of oil or natural gas.  It is also self-evident that the price of oil in commodity exchange markets fluctuates following divers macroeconomic and global geopolitical forces.
In this paper, we use the mean reverting regime switching L\'{e}vy processes to model the oil price. 
 Oil prices also display a great deal of seasonality, jumps, and spikes due to various supply disruptions and political turmoils in oil-rich countries. We use regime-switching jump diffusions to capture all those effects. Thus our pricing model closely captures the instability of oil markets. 
The paper is organized as follows. In the next section, we
formulate the problem under consideration. In Section 3, we prove the existence of a Nash equilibrium. And in section 4, we construct a finite difference approximation scheme and prove its convergence to the value functions. Finally, in section 5, we give a numerical example.

\section{Problem formulation}
Consider a multinational oil company who enters into a Production Sharing Agreement with the government of an oil-rich country with expiration $0<T<\infty$.  Both parties will share the profits from the sales of the extracted oil on world markets following a simple rule where the company takes $100\theta$ percent and the government takes $100(1-\theta)$ percent of the profits, for some $\theta\in (0,1)$ . We assume that  the market value of a barrel of oil at time $t$ is $S_t=e^{X_t}$. In fact, we assume that the oil price follows an exponential L\'{e}vy model, these models are natural extensions of the celebrated Black-Scholes model. Given that oil prices are very sensitive global macroeconomic and geopolitical shocks, we model $X_t$ as a  mean reverting regime switching L\'{e}vy process with two states. Let $\al(t)\in{\cal M}=\{1,2\}$ be a finite state Markov chain that captures the state of the oil market: $\al(t)= 1$ indicates the bull market at time $t$ and $\al(t) = 2$ represents a bear market at  time $t$. Let the matrix $Q$ be the generator of $\al(t)$.   
Let  $(\eta_t)_t$ be a L\'{e}vy process and let $N$ be the Poisson random measure of $(\eta_t)_t$. Let $\nu$ be the L\'{e}vy measure of $(\eta_t)_t$.
The differential form of $N$ is denoted by $N(dt,dz)$, we define the  differential $\bar{N}(dt,dz)$ as follows: $\bar{N}(dt,dz)= N(dt,dz)-\nu(dz)dt$ if  $ |z|<1$ and $\bar{N}(dt,dz)= N(dt,dz)$ if $|z|\geq 1.$
We assume that  the L\'{e}vy measure $\nu$ has finite intensity, $ 
\disp \Gamma := \int_{\rr} \nu(dz)<\infty.$
Let $ K<\infty$ be the total size of the oil field at the beginning of the lease, and let $Y(t)$ be the size of the remaining reserve of the oil field by time $t$, obviously $Y(t)\in [0,K]$. We model the evolution of the profit sharing agreement as a differential game where the two players are the oil company and the government. The state variables of our differential game are $X(t)\in \rr$ and $Y(t)\in [0,K]$, and the state space is $\rr\times[0,K]$, because the oil price $S_t$ is fully determined by its logarithm $X_t=\log(S_t)\in \rr$.
We will refer the oil company as Player 1 and the government as Player 2.
 We  assume that  the processes $ X(t)$ and $ Y(t)$ follow the   dynamics
\beq{model}
 \left \{ \begin{array}{ll}
 \mathrm{d}X(t)&\disp =\kappa \big(\mu(\al(t))- X(t)\big)\mathrm{d}t\\
&+\sigma(\alpha(t))\mathrm{d}W(t)  \\
 &+\int_{\rr}\gamma(\al(t))z\bar{N}(dt,dz),\\
 dY(t)&= -u_1(t) dt, \\
X(s)&=x, \,\,\,\, Y(s) = y \geq 0,\,\, 0\leq s\leq t\leq T, \end{array} \right.
\eeq   
where $ u_1(t)\in U_1=[0,\bar{u}_1]$ is the extraction rate chosen by the company and $ u_2(t) \in U_2= [0,\bar{u}_2]$ is the tax rate chosen by the government. The constant  $\bar{u}_1$ represents the maximum extraction rate and $\bar{u}_2$ is the maximum tax rate. The processes $ u_1(t)$ and $ u_2(t)$ are control variables, and $ W(t)$ is the Wiener process defined on a probability space $ (\Omega,{\cal F },P)$. Moreover, we assume that $W(t)$, $\al(t)$ and $\eta_t$ are independent. The  parameter $\exp(\mu(\cdot))$ represents the equilibrium price of oil and $\kappa$ represents the coefficient of mean reversion. For each state $i\in \{1,2\}$ of the oil market, we assume that the corresponding equilibrium price $\exp(\mu(i))$ is known. Similarly $\sigma(\cdot)$ represents the volatility and $\gamma(\cdot)$ represents and the intensity of the jump diffusion. For each state $i\in\{1,2\}$ of the oil market, we assume that $\sigma(i),\, \gamma(i)$ are known nonzero constants. As a matter of fact, $\gamma(i)$ captures the frequencies and jump sizes of the oil price.
\begin{defn}
The extraction and taxation rates  $u_1(\cdot)$ and $u_2(\cdot)$, taking values on intervals $[0,\bar{u}_1]$ and $[0,\bar{u}_2]$ respectively, are called admissible controls with respect to the initial data $(s,x,y, i)\in {\cal D}:= [0,T]\times\rr\times[0,K]\times{\cal M}$ if:
 \begin{itemize}
\item Equation (\ref{model}) has a unique  solution  with $X(s) = x$,  $Y(s) =y$, $\al(s)= i $, and $X(t)\in\rr, \,\, Y(t)\in[0,K]$ for all $t\in [0,T]$.
\item The processes $u_1(\cdot)$ and $u_2(\cdot)$ are $\{{\cal F}_t\}_{t\geq 0}$-adapted where $ {\cal F}_t = \sigma\{\alpha(s),W(s), \eta_s; s\leq t\}$.
\end{itemize}
We use ${\cal U}_j ={\cal U}_j(s,x,y,i)$ to denote the set of admissible controls taking values in $U_j=[0,\bar{u}_j]$ such that $X(s)=x$, $Y(s) =y$, $\al(s)=i)$,  for each $j\in\{1,2\}$. 
\end{defn}
Let $C(t,Y(t),u_1(t))$  be  the extraction cost function per unit of time $t$. A typical example of extraction cost function is 
$ C(t,y,u) := a+ m u(c-by),
$
where $a>0$ can be seen as the initial cost  of setting up the oil field and $m$,  $b$, and $c$ are constants such that $m>0$ and $b,c\geq 0$.
%
The total profit rate for operating the mine is
\bea
 P(t,X(t),Y(t),u_1(t))=e^{X(t)} u_1(t) -\\
 C(t,Y(t), u_1(t)).
\eea
 The total income tax the government levies on the oil company is $ u_2(t) \theta P(t,X(t), Y(t), u_1(t))$.
 The post-tax profit rate of the company is 
\bea
\ad L_1(t,X(t),Y(t),u_1(t), u_2(t))\\
\ad= \theta P(t,X(t),Y(t), u_1(t) )(1- u_2(t)), 
\eea
and the government profit rate function is
\bea
\ad L_2(t,X(t),Y(t),u_1(t), u_2(t))\\
\ad= (1-\theta) P(t,X(t),Y(t),u_1(t)) + \\
\ad u_2(t)\theta P(t,X(t),Y(t),u_1(t)).
\eea
We assume that at the end of the lease there are no extraction  revenues, therefore  the profit rate for the oil company could be zero or equal to the cost of closing the mine. We will denote by
 $\Phi_1(X(T),Y(T)) $ the terminal profit rate for the oil company. In most cases, $\Phi_1(x,y):=0$.
However, the terminal profit rate of the government is the market value of the remaining reserve,
$\Phi_2(X(T),Y(T)) =Y(T)(\exp(X(T))-m)$, where $m$ is the cost of extracting one barrel.
In sum, we will generally assume that the running profit rate and terminal profit rate functions $L_i, \Phi_i$, $i=1,2$ are Lipschitz continuous on bounded sets.
Given a discount rate $r> 0$, the payoff functional of Player $i=1,2$ is 
\bea
\ad J_i(s,x,y,\iota;u_1,u_2)\nonumber\\
\ad= E\Bigg[\int_s^Te^{-r(t-s)}   L_i(t,X(t),Y(t),u_1(t), u_2(t) ) dt  \\
\ad + e^{-r(T-s)}\Phi_i(X(T),Y(T))\bigg{|} X(s) = x, 
 Y(s) = y, \\
 \ad \al(s) = \iota\Bigg].
\eea
Each player wants to maximize its own payoff.
The company will try to maximize its payoff by adjusting the extraction rate $ u_1(\cdot)$, while the government will maximize its payoff by changing the tax rate $ u_2(\cdot)$. We, therefore, model this interaction as a noncooperative differential game.   Our goal is to find a noncooperative {  Nash equilibrium} $ (u_1^*, u_2^*)$  such that  
\bea
\ad  J_1(s,x,y,\iota;u_1^*,u_2^*)\geq J_1(s,x,y,\iota;u_1,u_2^*), \\
\ad\hbox{ for all } u_1(\cdot) \in {\cal U}_1(s,x,y,\iota),\\
\ad  J_2(s,x,y,\iota;u_1^*,u_2^*)\geq J_2(s,x,y,\iota;u_1^*,u_2), \\
\ad \hbox{ for all }u_2(\cdot) \in {\cal U}_2(s,x,y,\iota).
\eea
 In the next section we will prove the existence of a Nash equilibrium.
\section{Nash Equilibrium}

\begin{defn}
Let $(u_1^*, u_2^*)$ be  a Nash equilibrium of our differential game, the functions 
\beq{val1}
 V_1(s,x,y,\iota) = \sup_{u_1\in{\cal  U}_1}J_1(s,x,y,\iota; u_1,u_2^*) \nonumber \\
  V_2(s,x,y,\iota) = \sup_{u_2\in{\cal  U}_2}J_2(s,x,y,\iota; u_1^*,u_2) \nonumber
\eeq
are called value functions of Player 1 and Player 2 respectively.
\end{defn}
  In order  to find the optimal strategies $ u_1^*$ and $ u_2^*$ of a Nash equilibrium we first have to derive the value functions $ V_1$ and $ V_2$ of the differential game. 
  Formally the value functions $ V_1$ and $ V_2$ should satisfy  the following Hamilton-Jacobi-Isaacs equations.
  Assuming that we have a Nash equilibrium $ (u_1^*, u_2^*)$ let us define the corresponding Hamiltonians:
\beq{Hamilton1}
&& H_1\bigg(s,x,y,\iota,V,\frac{\partial V}{\partial x},\frac{\partial V}{\partial y}, \frac{\partial^2 V}{\partial x^2}\bigg)     \nonumber\\
&&=rV - \sup_{u_1\in U_1} \Bigg{(}\frac{1}{2}\sigma^2(\iota)\frac{\partial^2V}{\partial x^2}  
+  \kappa \big(\mu(\iota)- x\big)\frac{\partial V}{\partial x}\nonumber\\
 &&-  u_1\frac{\partial V}{\partial y}   
 + \nonumber  \int_{\rr}\bigg(V(s,x+\gamma(\iota)z,y,\iota) -V(s,x,y,\iota) \nonumber \\
&&  -{\bf 1}_{\{|z|<1\}}(z)\frac{\partial V}{\partial x} \gamma(\iota)z\bigg)\nu(dz)  \nonumber \\
&&+ L_1(s,x,y,u_1,u_2^*)+  QV(s,x,y,\cdot)(\iota)\Bigg{)},
\eeq
and	 
\beq{Hamilton2}
&& H_2\bigg(s,x,y,\iota,V,\frac{\partial V}{\partial x},\frac{\partial V}{\partial y}, \frac{\partial^2 V}{\partial x^2}\bigg)     \nonumber\\
&&=rV  - \sup_{u_2\in U_2} \Bigg{(}\frac{1}{2}\sigma^2(\iota)\frac{\partial^2V}{\partial x^2}  
+  \kappa \big(\mu(\iota)- x\big) \frac{\partial V}{\partial x} \nonumber \\
&&- u_1^*\frac{\partial V}{\partial y}   + \nonumber  \int_{\rr}\bigg(V(s,x+\gamma(\iota)z,y,\iota) -V(s,x,y,\iota) \nonumber \\
&& -{\bf 1}_{\{|z|<1\}}(z)\frac{\partial V}{\partial x}\cdot \gamma(\iota)z\bigg)\nu(dz) \nonumber \\
&&+ L_2(s,x,y,u_1^*,u_2)+  QV(s,x,y,\cdot)(\iota)\Bigg{)},
\eeq	 
with $\disp  QV(s,x,y,\cdot)(\iota) = \sum_{j\neq \iota}q_{\iota j}(V(s,x,y,j) -V(s,x,y,\iota))$. 
  The corresponding Hamilton Jacobi Isaacs equations of this noncooperative game are
\beq{isaac1}
 \left \{ \begin{array}{ll} \frac{\partial V_1}{\partial s}  =   H_1\bigg(s,x,y,\iota,V_1,\frac{\partial V_1}{\partial x},\frac{\partial V_1}{\partial y}, \frac{\partial^2 V_1}{\partial x^2}\bigg),\\
 \frac{\partial V_2}{\partial s} =   H_2\bigg(s,x,y,\iota,V_2,\frac{\partial V_2}{\partial x},\frac{\partial V_2}{\partial y}, \frac{\partial^2 V_2}{\partial x^2}\bigg),\\
V_1(T,x,y,\iota) = \Phi_1(x,y)   \\
V_2(T,x,y,\iota) = \Phi_2(x,y).
\end{array}\right.
\eeq
The strategy for solving this differential game is to first find the solutions of the Isaacs equation (\ref{isaac1}) and then derive the optimal extraction and taxation policies from  the Nash equilibrium.  
The next result gives the road map we will use to find a Nash equilibrium if we already have the value functions.
\begin{thm}\label{NashEquilibrium}
Assume that there exists $ (u_1^*, u_2^*)\in {\cal U}_1\times {\cal  U}_2$ such that 
the nonlinear Hamilton-Jacobi-Isaacs equations (\ref{isaac1})  have classical  solutions $ V_i(s,x,y,\iota)$, $i=1,2$ with
 \beq{optControl1} 
  \disp  u_1^*(s) = \arg\max\bigg(   -  u_1\frac{\partial V_1(s,x,y,\iota)}{\partial y}    \nonumber\\
  +  L_1(s,x,y,u_1,u^*_2)  \bigg), \quad s\in[0,T]\nonumber 
\eeq
and
 \beq{optControl2} 
  \disp  u_2^* (s)= \arg\max\bigg(  \disp  L_2(s,x,y,u_1^*,u_2) \bigg),\, s\in[0,T].  \nonumber 
\eeq
Then the pair $  (u^*_1, u^*_2)$ is a Nash equilibrium  and 
$  J_i(s,x,y,\iota;u_1^*,u_2^*)  = V_i(s,x,y,\iota)$, $i=1,2$.
\end{thm}
\paragraph{Proof.}
The proof relies on the fact that this problem can be  uncoupled and solved as  an optimal control problem.  In fact, if we  replace the control process  $u_2(\cdot)$ by $ u^*_2(\cdot)$ in (\ref{model}) then, the differential game problem becomes an optimal control problem with the only control variable $u_1(\cdot)$.  The HJB equation of this new control problem is 
\beq{hjb1}
 \left \{ \begin{array}{ll} \frac{\partial W_1}{\partial s}  =   H_1\bigg(s,x,y,\iota,W_1,\frac{\partial W_1}{\partial x},\frac{\partial W_1}{\partial y}, \frac{\partial^2 W_1}{\partial x^2}\bigg)&\\
W_1(T,x,y,\al(T)) = \Phi_1(x,y).
\end{array}\right.
\eeq
Following the assumptions of this theorem, it is clear that the HJB equation (\ref{hjb1}) has a solution $V_1$ and the optimal policy of this new control problem is $u_1^*$. Therefore $u_1^*$ is in equilibrium with $u_2^*$ and $V_1$ is the value function of  Player 1. A similar argument can be used to show that $u_2^*$ is in equilibrium with $u_1^*$ and that $V_2$ is the value function of  Player 2.
\endproof
Given that $L_i, \Phi_i, i=1,2$ are Lipschitz continuous, using standard methods from control theory it can be shown that the values functions $V_1$ and $V_2$ are the unique viscosity solutions of the Isaacs equations (\ref{isaac1}).
 The uniqueness of the viscosity solutions are obtained as in   Barles and Imbert (2008)  by applying nonlocal extensions of the Jensen-Ishii Lemma. For more on the derivation of the maximum principle for nonlocal operators, one can also refer to Biswas et al. (2010). 
\section{Numerical Approximation}
In this section, we construct a finite difference scheme and
show that it converges to the unique viscosity solutions of the Isaacs equation (\ref{isaac1}).  We will use the following notations; we set $u=(u_1,u_2)$, $ u^*=(u_1^*, u_2^*)$, $ U=U_1\times U_2$, and 
  \bea
V(s,x,y,i)= \left ( \begin{array}{ll}  V_1(s,x,y,i)\\
V_2(s,x,y,i) \end{array}\right),
\eea
  \bea
\ad \disp H(s,x,y,i,V,\frac{\partial V}{\partial x},\frac{\partial V}{\partial y}, \frac{\partial^2 V}{\partial x^2}) \\
\ad= \left ( \begin{array}{ll}\disp    H_1(s,x,y,i,V_1,\frac{\partial V_1}{\partial x},\frac{\partial V_1}{\partial y}, \frac{\partial^2 V_1}{\partial x^2}) \\
\disp  H_2(s,x,y,i,V_2,\frac{\partial V_2}{\partial x},\frac{\partial V_2}{\partial y}, \frac{\partial^2 V_2}{\partial x^2})  \end{array}\right), 
\eea
  \bea
\Phi(x,y)= \left ( \begin{array}{ll} \Phi_1(x,y)\\
\Phi_2(x,y) \end{array}\right). 
\eea
The Isaacs equation (\ref{isaac1}) can be rewritten as follows
\beq{isaac2}
 \left \{ \begin{array}{ll}\disp \frac{\partial V}{\partial s}  =   H\bigg(s,x,y,i,V,\frac{\partial V}{\partial x},\frac{\partial V}{\partial y}, \frac{\partial^2 V}{\partial x^2}\bigg), \\
V(T,x,y,i) = \Phi(x,y).
\end{array}\right.
\eeq
Let $k\in(0,1) $ be the step size with respect to $s$, and $h\in (0,1)$ be the step size with respect to  $x$ and $y$,  we will use the standard finite difference
operators  $\Delta_s$, $\Delta_x$, $\Delta_{xx}$ and $\Delta_y$.
Let ${ \rm I}f$ denote the integral part of the Hamiltonians $H_1$ and $H_2$. We will approximate ${\rm I}f$ using the Simpson's quadrature.
Using the fact the L\'{e}vy measure is finite  $\disp \Gamma = \int_{\rr} \nu(dz)<\infty$, we have 
\bea
\ad{\rm I}f(s,x,y,i)= \int_{\rr}f(s,x+\gamma(i)z,y,i)\nu(dz)\\
\ad-\frac{\partial f(s,x,y,i)}{\partial x}\int_{-1}^1 \gamma(i)z\nu(dz)  -f(s,x,y,i)\Gamma. 
\eea
We use the Simpson's quadrature to approximate the integral part of the Hamiltonians.
Let $\xi\in(0,1)$ be the step size of the Simpson's quadrature, the corresponding approximation of the the integral part is
\bea
\ad{\rm I}_{\xi}f(s,x,y,i) = \sum_{j=0}^{N_\xi}c_jf(s,x+\gamma(i)z_j,y,i)-\\
\ad\frac{\partial f(s,x,y,i)}{\partial x}\sum_{j=0}^{M_\xi} d_j\gamma(i)z_j-f(s,x,y,i)\Gamma, 
\eea
where the $(c_j)_{0\leq j\leq N_\xi}$ and $(d_j)_{0\leq j\leq M_\xi}$ are the corresponding sequences of the coefficients of the Simpson's quadrature and $z_j\in[-1,1]$, $j\in\{1,...,M_\xi\}$ are grid points of the interval $[-1,1]$. In fact, $\disp  \lim_{N_\xi \rightarrow \infty}  \sum_{j=0}^{N_\xi}c_j =\Gamma$ and $\disp \lim_{M_\xi \rightarrow \infty}\sum_{j=0}^{M_\xi}d_j=\int_{-1}^1\nu(dz) $.
The corresponding discrete versions of the Hamiltonians $H_1, H_2$ are defined as follows
\beq{disc-Ham1}
&&  H_{1, u_2^*}^{h,k,\xi}V_1(s,x,y,i) \nonumber \\
 &=&rV_1(s,x,y,i) -  \sup_{u_1\in U_1} \Bigg(\frac{1}{2}\sigma^2(i)\Delta_{xx}V_1(s,x,y,i)\nonumber\\
 &&+{\rm I}_{\xi}V_1(s,x,y,i)+ \kappa( \mu(i)-x)\Delta_x V_1(s,x,y,i)  \nonumber\\
&&
 - u_1\Delta_yV_1(s,x,y,i)+L_1(s,x,y,u_1,u_2^*) \nonumber\\
&&
  + QV_1(s,x,y,\cdot)(i)\Bigg)  \nonumber
\eeq
and
\beq{disc-Ham2}
&&  H_{2, u_1^*}^{h,k,\xi}V_2(s,x,y,i) \nonumber \\
 &=& rV_2(s,x,y,i) -  \sup_{u_2\in U_2} \Bigg(\frac{1}{2}\sigma^2(i)\Delta_{xx}V_2(s,x,y,i)  \nonumber \\
 &&+{\rm I}_{\xi}V_2(s,x,y,i)+ \kappa( \mu(i)-x)\Delta_x V_2(s,x,y,i) \nonumber\\
&&- u_1^*\Delta_yV_2(s,x,y,i)+L_1(s,x,y,u_1^*,u_2)
 \nonumber\\
&& + QV_2(s,x,y,\cdot)(i)\Bigg). \nonumber
\eeq
Therefore the discrete version of (\ref{isaac2}) is 
\beq{discrete}
 \left \{ \begin{array}{ll}\disp 0 = \Delta_s V  - H_{u_1^*,u_2^*}^{h,k,\xi}V(s,x,y,i),& \\
V(T,x,y,i) = \Phi(x,y), 
\end{array}\right.
\eeq
with
\bea
H_{u_1^*,u_2^*}^{j,k,\xi}V(s,x,y,i) =\left ( \begin{array}{ll} H_{1, u_2^*}^{h,k,\xi}V_1(s,x,y,i)\\
 H_{2,u_1^*}^{h,k,\xi}V_2(s,x,y,i) 
\end{array}\right).
\eea
 We have the following crucial Lemma.  

\begin{lem}\label{DiscreteLemma}
Let $h\in(0,1) $ be small enough, for each $k, \xi \in (0,1)$, there exists a  unique bounded function $V_{h,k, \xi}$  defined on ${\cal D}$  that solves equation (\ref{discrete}).
\end{lem}
\paragraph{Proof}
Any solution $V$ of (\ref{discrete}) should also satisfy $V =V + \epsilon(\Delta_s - H_{u_1^*,u_2^*}^{h,k,\xi}V(s,x,y,i))$ for $\epsilon>0$ together with terminal condition $V(T,x,y,i) = \Phi(x,y)$.
We define the operator ${\cal F}_\xi$  on bounded functions on ${\cal D}$ as follows
\beq{DiscrOpe}
& & {\cal F}_\xi(V)( s,x,y,i; h,k)\nonumber \\
&=& V + \epsilon(\Delta_s - H_{u_1^*,u_2^*}^{h,k}V(s,x,y,i))\nonumber\\
&=&  \sup_{u\in U}\Bigg(\frac{\epsilon }{k}V(s+k,x,y,i) +  \epsilon  PV(s,x,y,\cdot)(i)  \nonumber \\
&&+ b(i)V(s,x-h,y,i) +c_{u,u^*}(i)V(s,x,y,i)  \nonumber\\
&& 
+ \epsilon L(s,x,y,u_1,u_2)
+ a_{u,u^*}(i)V(s,x+h,y,i)  \nonumber\\
&&+\epsilon \sum_{n=0}^{N_\xi} c_n V(s,x+\gamma(i)z_n,y,i)\nonumber
\Bigg),\\
&&{\cal F}_\xi(V)(T,x,y,i;h,k) = \Phi(x,y),
\eeq 
where the coefficients $q_{in}$ are coefficients of the generator $ Q$ and the quantities   $a_{u,u^*}(i)$, $b(i)$ ,  $c_{u,u^*}(i)$ and $ PV(s,x,y,\cdot)(i) $  are defined as follows 
\bea
 e(i;v) \ad= 1-\epsilon\bigg[r+\frac{1}{ k} +  \frac{\sigma^2(i)}{h^2} +  \frac{1}{h}\bigg(\kappa( \mu(i)-x) \\
\ad- \sum_{n=0}^{M_\xi} d_n\gamma(i)z_n
 - v\bigg) + \Gamma  +\sum_{n\ne i} q_{in}\bigg] ,
 \eea
 \bea
 g(i;v) \ad=  \frac{\sigma^2(i)}{2h^2} +\frac{1}{h}\bigg(\kappa( \mu(i)-x)\\
 \ad\hspace{0.5in}- \sum_{n=0}^{M_\xi} d_n\gamma(i)z_n
 - v\bigg),
 \eea
 \bea
\ad c_{u,u^*}(i) =\left(\begin{array}{ll}e(i;u_1)  & 0 \\
0 &  e(i;u_1^*)
\end{array}\right),
\eea
\bea
\ad  PV(s,x,y,\cdot)(i) =\left(\begin{array}{ll}\disp  \sum_{n\ne i} q_{in} V(s,x,y,n) \quad 0\\
\disp 0\quad  \sum_{n\ne i} q_{in}V(s,x,y,n)
\end{array}\right), 
\eea
\bea
a_{u,u^*}(i)
\ad =\epsilon\left(\begin{array}{ll} g(i,u_1) \quad 0\\
0 \quad   g(i,u_1^*)\end{array}\right),
\eea
\bea
\ad b(i)  =\frac{\sigma^2(i)\epsilon}{2h^2}I, \,\ \,\,\, I=\left(\begin{array}{ll}1 & 0\\
0&1
\end{array}\right).
 \eea
Note that equation (\ref{discrete}) is equivalent to $ V(s,x,y,i) = {\cal F}_\xi(V)(s,x,y,i;h,k)$, it suffices to show the operator ${\cal F}_\xi$ has a fixed point.  As a matter of fact, for $h$ small enough, it is clear that  $g(i,v)>0$ for all  $i\in {\cal M}$ and  $v\in U_1$, thus  both diagonal coefficients of $ a_{u,u^*}(i)$ are positive. Moreover, we can choose  $\epsilon $ such that $e(v,i)>0,$  for all  $v\in U_1, i\in {\cal M}$ and 
\bea
\ad  a_{u,u^*}(i) +b(i)+ c_{u,u^*}(i) + \epsilon \bigg(\frac{1}{k}  +  \sum_{n=0}^{N_\xi} c_n\\
\ad+ \sum_{n\ne i} q_{in} \bigg)I   = \bigg(1-\epsilon r+\epsilon\big( \sum_{n=0}^{N_\xi} c_n-\Gamma\big)\bigg)I,\\
\hbox{with }\ad   0\leq 1-\epsilon r+\epsilon\big( \sum_{n=0}^{N_\xi} c_n-\Gamma\big)\leq \delta<1,
\eea 
where $\delta$ is  a constant.
Using the fact that the difference of two suprema is less than the supremum of the difference.  If we have two bounded functions $V, W $  defined on $ {\cal D}$, it is clear that 
\bea
\ad |{\cal F}_\xi(V)(s,x,y,i;h,k) - {\cal F}_\xi(W)(s,x,y,i;h,k)|  \\
\ad \leq \delta \sup_{{\cal D}}|V-W|.
\eea
Therefore,  the map ${\cal F}_\xi$ is a contraction on the space of bounded functions on $ {\cal D}$, using the Banach's Fixed Point Theorem we conclude the proof of the lemma.
\endproof
\begin{rem}
It is clear from \lemref{DiscreteLemma} that the numerical scheme obtained from (\ref{discrete}) is stable  since the solution of the scheme is bounded independently of the step sizes $h,k, \xi\in(0,1)$ and obviously consistent because as the step sizes $h,k,\xi$ go to zero the finite difference operators converge to the actual partial differential operators. We have the following convergence theorem.
\end{rem}
\begin{thm}\label{Convergence}
 Let $V_{h,k,\xi}$ be the solution of the discrete scheme obtained in \lemref{DiscreteLemma}. Then as  $(h,k,\xi)\rightarrow 0$ the sequence $V_{h,k,\xi}$ converges locally uniformly on $ {\cal D}$ to the unique viscosity solution $V$ of (\ref{isaac2}).
\end{thm}
This result  is the standard method for approximating  viscosity solutions, for more one can refer to Barles and Souganidis (1991).

\section{Applications}
Consider an oil company with a 10 years lease to extract oil from an oil field with an known  capacity of K=10 billion barrels. We assume that the profit sharing agreement between the oil company  and the government is such that the oil company takes 40\%  of  profits and the government takes 60\%, so $\theta = 0.4$.  The oil equilibrium price when the market is up is $\mu(1)=50$ and when the market is down is $\mu(2) = 35$. The mean reversion coefficient is $\kappa = 0.01$,
the volatility  when the market is up is $\sigma(1) = 0.1$ and when the market is down $ \sigma(2) =0.3$.
 And the jump intensity is $ \lambda(1) =0.01$ when the market is up and $ \lambda(2)=0.15$ when the market is down. The generator of the Markov chain is $\disp  Q= \left ( \begin{array}{ll}   -0.005 & 0.005\\
 0.002& -0.002\end{array}\right)$. We assume that $
L_1(t,x,y,u_1,u_2) = 0.4(e^xu_1-(10+15 u_1))(1-u_2)$ ,  $\Phi_1(x,y)=0$, 
$
L_2(t,x,y,u_1,u_2) = 0.6(e^xu_1-(10+15 u_1))+ 0.4u_2(e^xu_1-(10+15 u_1)),$ and  $\Phi_2(x,y) = y(e^x-15)$.
Moreover, we assume that the extraction $ u_1(\cdot) \in [0,50000]$ and the top tax rate is { 30\% }  so $ u_2(\cdot)\in [0,0.3]$.  Keep in mind that, because the payoff rates are linear functions of each control variable $ u_1(\cdot)$ and $u_2(\cdot)$, therefore using \thmref{NashEquilibrium}, the optimal strategies $ u_1^*$ and $u_2^*$  are obtained by looking at the signs
of the following functionals $
\ F(t,x,y,u_2,i) =  -\frac{\partial V_1(t,x,y,i)}{\partial y} + 0.4 (e^x - 15))(1- u_2(t))$, and $
 G(t,x,y,u_1,i)  =  0.4 (e^xu_1(t) - (10 +15u_1(t)) ).
$
The optimal strategies will only be attained at the endpoints of the intervals $ U_1=[0,50000] $ and $ U_2= [0,0.3]$, we have
\bea 
 u_1^*(t) = \left \{ \begin{array}{ll} 0 & \hbox{ if    }\quad  F(t,x,y,u_2^*(t),i)   \leq 0\\
50000 &\hbox{ if    }  \quad  F(t,x,y,u_2^*(t),i)  >0 ,\end{array}\right.
\eea
and 
\bea 
 u_2^*(t) = \left \{ \begin{array}{ll} 0 & \hbox{ if } \quad  G(t,x,y,u_1^*(t),i) \leq 0\\
0.3  &\hbox{ if }  \quad G(t,x,y,u_1^*(t),i)>0 .\end{array}\right.
\eea

In the next two figures we have the plots of the functions $F$ and $G$ when the market is up and when the market is down. For instance, in Figure 1 we have plots of the $F$ and $G$ when the market is bullish, the regions above the curves represent the domains where signs of $F$ and $G$ are positive, so in those regions it is always optimal to extract at full capacity  or to tax at the maximal rate. And the regions below the curves represent the domains where the signs of $F$ and $G$ are negative, thus in those regions it is optimal not to extract at all and not to tax the oil company.  In Figure 2, similar plots are given when the market is bearish.

\begin{figure}
\vspace{-2cm}
\includegraphics[ height=10cm,width=11cm]{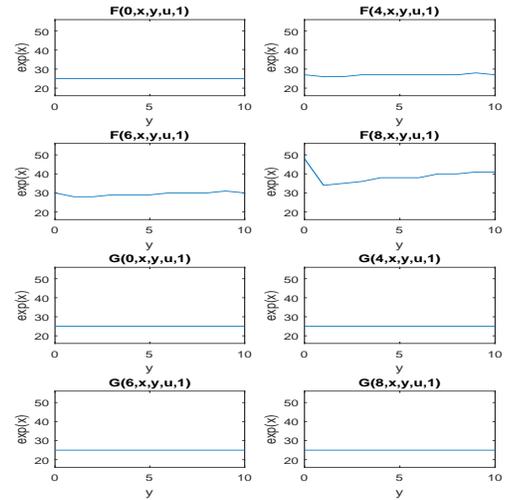}
\caption{This graph represents the optimal regions when the market is bullish. }
\end{figure}

\begin{figure}
\vspace{-2cm}
\includegraphics[ height=10cm,width=11cm]{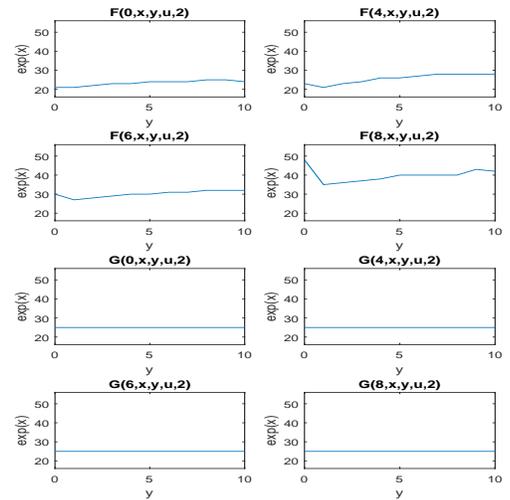}
\caption{This graph represents the optimal regions when the market is bearish }
\end{figure}




\begin{thebibliography}{17}



\bibitem{Bar2}
G. Barles and C. Imbert, Second-order  elliptic integro-differential  equations: viscosity solutions' theory revisited, {\it Ann. Inst. Henri Poincar\'{e} Anal. Non Lin\'{e}aire}, {\bf 25}, 3, (2008), pp. 567-585.


\bibitem{Bar}
G. Barles and P.E. Souganidis,
Convergence of approximation schemes for fully nonlinear second order
equations, {\it Asymptot. Anal.} {\bf 4}, (1991), 271-283 


\bibitem{Biswas}
I. H. Biswas, E. R. Jakobsen and K. H. Karlsen, {Viscosity  solutions for a system of integro-PDE and connections to optimal switching and control of jump-diffusion processes}, {\it Applied Mathematics and Optimization}, {\bf 62}, (2010), pp. 47-80.






\bibitem{Hotelling}
H. Hotelling,  The economics of exhaustible resources, {\it Journal of Political Economy},  {\bf 39}, 2, (1931), pp. 137-175.  



























\end{thebibliography}
\end{document}